\documentclass[11pt]{article}
\usepackage{amsmath,amssymb,amsbsy,amsfonts,amsthm,latexsym,amsopn,amstext,
                                    amsxtra,euscript,amscd}
\begin{document}




\newfont{\teneufm}{eufm10}
\newfont{\seveneufm}{eufm7}
\newfont{\fiveeufm}{eufm5}
%
%
\newfam\eufmfam
                    \textfont\eufmfam=\teneufm \scriptfont\eufmfam=\seveneufm
                    \scriptscriptfont\eufmfam=\fiveeufm

%
%
\def\frak#1{{\fam\eufmfam\relax#1}}
%


\def\bbbr{{\rm I\!R}} 
\def\bbbc{{\rm I\!C}} 
\def\bbbm{{\rm I\!M}}
\def\bbbn{{\rm I\!N}} 
\def\bbbf{{\rm I\!F}}
\def\bbbh{{\rm I\!H}}
\def\bbbk{{\rm I\!K}}
\def\bbbl{{\rm I\!L}}
\def\bbbp{{\rm I\!P}}
\newcommand{\lcm}{{\rm lcm}}
\def\bbbone{{\mathchoice {\rm 1\mskip-4mu l} {\rm 1\mskip-4mu l}
{\rm 1\mskip-4.5mu l} {\rm 1\mskip-5mu l}}}
\def\bbbc{{\mathchoice {\setbox0=\hbox{$\displaystyle\rm C$}\hbox{\hbox
to0pt{\kern0.4\wd0\vrule height0.9\ht0\hss}\box0}}
{\setbox0=\hbox{$\textstyle\rm C$}\hbox{\hbox
to0pt{\kern0.4\wd0\vrule height0.9\ht0\hss}\box0}}
{\setbox0=\hbox{$\scriptstyle\rm C$}\hbox{\hbox
to0pt{\kern0.4\wd0\vrule height0.9\ht0\hss}\box0}}
{\setbox0=\hbox{$\scriptscriptstyle\rm C$}\hbox{\hbox
to0pt{\kern0.4\wd0\vrule height0.9\ht0\hss}\box0}}}}
\def\bbbq{{\mathchoice {\setbox0=\hbox{$\displaystyle\rm
Q$}\hbox{\raise
0.15\ht0\hbox to0pt{\kern0.4\wd0\vrule height0.8\ht0\hss}\box0}}
{\setbox0=\hbox{$\textstyle\rm Q$}\hbox{\raise
0.15\ht0\hbox to0pt{\kern0.4\wd0\vrule height0.8\ht0\hss}\box0}}
{\setbox0=\hbox{$\scriptstyle\rm Q$}\hbox{\raise
0.15\ht0\hbox to0pt{\kern0.4\wd0\vrule height0.7\ht0\hss}\box0}}
{\setbox0=\hbox{$\scriptscriptstyle\rm Q$}\hbox{\raise
0.15\ht0\hbox to0pt{\kern0.4\wd0\vrule height0.7\ht0\hss}\box0}}}}
\def\bbbt{{\mathchoice {\setbox0=\hbox{$\displaystyle\rm
T$}\hbox{\hbox to0pt{\kern0.3\wd0\vrule height0.9\ht0\hss}\box0}}
{\setbox0=\hbox{$\textstyle\rm T$}\hbox{\hbox
to0pt{\kern0.3\wd0\vrule height0.9\ht0\hss}\box0}}
{\setbox0=\hbox{$\scriptstyle\rm T$}\hbox{\hbox
to0pt{\kern0.3\wd0\vrule height0.9\ht0\hss}\box0}}
{\setbox0=\hbox{$\scriptscriptstyle\rm T$}\hbox{\hbox
to0pt{\kern0.3\wd0\vrule height0.9\ht0\hss}\box0}}}}
\def\bbbs{{\mathchoice
{\setbox0=\hbox{$\displaystyle     \rm S$}\hbox{\raise0.5\ht0\hbox
to0pt{\kern0.35\wd0\vrule height0.45\ht0\hss}\hbox
to0pt{\kern0.55\wd0\vrule height0.5\ht0\hss}\box0}}
{\setbox0=\hbox{$\textstyle        \rm S$}\hbox{\raise0.5\ht0\hbox
to0pt{\kern0.35\wd0\vrule height0.45\ht0\hss}\hbox
to0pt{\kern0.55\wd0\vrule height0.5\ht0\hss}\box0}}
{\setbox0=\hbox{$\scriptstyle      \rm S$}\hbox{\raise0.5\ht0\hbox
to0pt{\kern0.35\wd0\vrule height0.45\ht0\hss}\raise0.05\ht0\hbox
to0pt{\kern0.5\wd0\vrule height0.45\ht0\hss}\box0}}
{\setbox0=\hbox{$\scriptscriptstyle\rm S$}\hbox{\raise0.5\ht0\hbox
to0pt{\kern0.4\wd0\vrule height0.45\ht0\hss}\raise0.05\ht0\hbox
to0pt{\kern0.55\wd0\vrule height0.45\ht0\hss}\box0}}}}
\def\bbbz{{\mathchoice {\hbox{$\sf\textstyle Z\kern-0.4em Z$}}
{\hbox{$\sf\textstyle Z\kern-0.4em Z$}}
{\hbox{$\sf\scriptstyle Z\kern-0.3em Z$}}
{\hbox{$\sf\scriptscriptstyle Z\kern-0.2em Z$}}}}
\def\ts{\thinspace}

\newtheorem{theorem}{Theorem}
\newtheorem{lemma}[theorem]{Lemma}
\newtheorem{claim}[theorem]{Claim}
\newtheorem{cor}[theorem]{Corollary}
\newtheorem{prop}[theorem]{Proposition}
\newtheorem{definition}{Definition}
\newtheorem{question}[theorem]{Open Question}
\newtheorem{remark}[theorem]{Remark}

\def\squareforqed{\hbox{\rlap{$\sqcap$}$\sqcup$}}
\def\qed{\ifmmode\squareforqed\else{\unskip\nobreak\hfil
\penalty50\hskip1em\null\nobreak\hfil\squareforqed
\parfillskip=0pt\finalhyphendemerits=0\endgraf}\fi}

\def\cA{{\mathcal A}}
\def\cB{{\mathcal B}}
\def\cC{{\mathcal C}}
\def\cD{{\mathcal D}}
\def\cE{{\mathcal E}}
\def\cF{{\mathcal F}}
\def\cG{{\mathcal G}}
\def\cH{{\mathcal H}}
\def\cI{{\mathcal I}}
\def\cJ{{\mathcal J}}
\def\cK{{\mathcal K}}
\def\cL{{\mathcal L}}
\def\cM{{\mathcal M}}
\def\cN{{\mathcal N}}
\def\cO{{\mathcal O}}
\def\cP{{\mathcal P}}
\def\cQ{{\mathcal Q}}
\def\cR{{\mathcal R}}
\def\cS{{\mathcal S}}
\def\cT{{\mathcal T}}
\def\cU{{\mathcal U}}
\def\cV{{\mathcal V}}
\def\cW{{\mathcal W}}
\def\cX{{\mathcal X}}
\def\cY{{\mathcal Y}}
\def\cZ{{\mathcal Z}}

\newcommand{\comm}[1]{\marginpar{%
\vskip-\baselineskip 
\raggedright\footnotesize
\itshape\hrule\smallskip#1\par\smallskip\hrule}}





\def\MOV{{\bf{MOV}}}

\hyphenation{re-pub-lished}

\def\ord{{\mathrm{ord}}}
\def\Nm{{\mathrm{Nm}}}
\renewcommand{\vec}[1]{\mathbf{#1}}

\def \F{{\bbbf}}
\def \L{{\bbbl}}
\def \K{{\bbbk}}
\def \Z{{\bbbz}}
\def \N{{\bbbn}}
\def \Q{{\bbbq}}
\def\E{{\mathbf E}}
\def\H{{\mathbf H}}
\def\G{{\mathcal G}}
\def\O{{\mathcal O}}
\def\cS{{\mathcal S}}
\def \R{{\bbbr}}
\def\Fp{\F_p}
\def \fp{\Fp^*}
\def\\{\cr}
\def\({\left(}
\def\){\right)}
\def\fl#1{\left\lfloor#1\right\rfloor}
\def\rf#1{\left\lceil#1\right\rceil}

\def\Zm{\Z_m}
\def\Zt{\Z_t}
\def\Zp{\Z_p}
\def\Um{\cU_m}
\def\Ut{\cU_t}
\def\Up{\cU_p}

\def\ep{{\mathbf{e}}_p}
\def\HH{\cH}

\def \Prob{{\mathrm {}}}

\def\LC{{\cL}_{C,\cF}(Q)}
\def\LCn{{\cL}_{C,\cF}(nG)}
\def\Tr{\mathrm{Tr}\,}

\def\taubar{\overline{\tau}}
\def\sigmabar{\overline{\sigma}}
\def\Fn{\F_{q^n}}
\def\En{\E(\Fn)}

\def\mand{\qquad \mbox{and} \qquad}

\newcommand{\ZZ}{{\mathbb Z}}
\newcommand{\PP}{{\mathbb P}^1}

\newcommand{\FF}{{\mathbb F}}
\newcommand{\QQ}{{\mathbb Q}}
\newcommand{\fq}{{\FF_q}}
\newcommand{\fqstar}{{\FF^*_q}}
\newcommand{\fqm}{{\FF_{q^m}}}
\newcommand{\fqn}{{\FF_{q^n}}}
\newcommand{\fqtwo}{{\FF_{q^2}}}
\newcommand{\fqnstar}{{\FF^*_{q^n}}}
\newcommand{\ffive}{{\FF_5}}
\newcommand{\fqnmi}{{\FF_{q^{nm_i}}}}
\newcommand{\fqs}{{\FF_{q^s}}}
\newcommand{\ftwo}{{\FF_2}}
\newcommand{\ftwon}{{\FF_{2^n}}}
\newcommand{\ztwo}{{\ZZ_2}}
\newcommand{\qtwo}{{\QQ_2}}
\newcommand{\fthreen}{{\FF_{3^n}}}
\newcommand{\fthree}{{\FF_3}}
\newcommand{\fthreem}{{\FF_{3^m}}}
\newcommand{\fthreeseven}{{\FF_{3^7}}}
\newcommand{\fthreer}{{\FF_{3^{2^r}}}}

\newcommand{\Disc}{\operatorname{Disc}}
\newcommand{\Res}{\operatorname{Res}}
\newcommand{\wt}{\mbox{wt}}
\newtheorem{conj}[theorem]{Conjecture}

\newcommand{\hlinepad}{\raisebox{0pt}[2.5ex]{\strut}}
\newcommand{\Hline}{\hline\hlinepad}
\newcommand{\HHline}{\hline\hline\hlinepad}


\title{Generalization of a Theorem of Carlitz}

\author{
{\sc Omran Ahmadi}  \\
{Claude Shannon Institute}\\
{University College Dublin} \\
{Dublin 4, Ireland} \\
{\tt omran.ahmadi@ucd.ie}}

\date{\today}

\maketitle

\begin{abstract}
We generalize Carlitz' result on the number of self reciprocal monic
irreducible polynomials over finite fields by showing that similar explicit
formula hold for the number of irreducible polynomials obtained by a fixed
quadratic transformation. Our main tools are a combinatorial argument and 
Hurwitz genus formula.
\end{abstract}

\section{Introduction}
Let $\fq$ denote the finite field with
 $q$ elements, where $q$ is a prime power, and let $\fq[x]$ denote the
polynomial ring over $\fq$. For
$f(x)$, a polynomial of degree $m$ over $\fq$ whose constant term
is nonzero, its \emph{reciprocal} is the polynomial $f^*(x)=x^m
f(1/x)$ of degree $m$ over $\fq$. A polynomial $f(x)$ is called
\emph{self-reciprocal} if $f^*(x)=f(x)$. The reciprocal of an
irreducible polynomial is also irreducible. The roots of the
reciprocal polynomial are the reciprocals of the roots of the
original polynomial, and hence, any self-reciprocal irreducible
monic ({\it{srim}}) polynomial of degree greater than one must have even
degree, say $2n$.

Self-reciprocal irreducible polynomials over finite fields have
been studied by many authors. Carlitz \cite{Carlitz} counted the
number of srim polynomials of degree $2n$ over a finite field for
every $n$. He showed the following.
\begin{theorem}\cite{Carlitz}
Let $SRIM(2n,q)$ denote the number of srim polynomials of degree
$2n$ over $\fq$. Then
\[
SRIM(2n,q) = \left\{
    \begin{array}{clcl}
    \frac{1}{2n}(q^n-1), & \;\;\;\mbox{if $q$ is odd and\;} n=2^m, \\\\
    \frac{1}{2n}\displaystyle\sum_{\stackrel{d\mid n,}{d\;\mbox{\tiny odd}}}\mu(d) q^{n/d}, &\;\;\;\mbox{otherwise}.\\
    \end{array}
    \right.
\]
\end{theorem}
In \cite{Cohen-1} and \cite{Meyn-1}, Cohen and Meyn, respectively, obtained the same
result by methods simpler than that was used by Carlitz in ~\cite{Carlitz}.

 In~\cite{AhVeg}, it has been shown that
 if $K$ is a field and $p(x)\in K[x]$ is a self-reciprocal
polynomial of degree $2n$, then for some $f(x)\in K[x]$ of degree
$n$, we have
\begin{equation}\label{eqdos}
p(x)=x^nf(x+x^{-1})=x^nf\left(\frac{x^2+1}{x}\right) \; .
\end{equation}

This implies that self reciproal irreducible monic polynomials
can be studied in the context of quadratic transformation of irreducible polynomials. 
A quadratic transformation of an irreducible polynomial $f(x)$ of degree 
$n$ over $\fq$ is the polynomial
\begin{equation}
p(x)=h(x)^nf\left(\frac{g(x)}{h(x)}\right),
\end{equation}
where $g(x)$ and $h(x)\in \fq[x]$ are coprime polynomials
with 
$$\max(\deg(g),\deg(h))=2.$$

Now it is natural to ask whether there exists a similar explicit formula for
the number of irreducible polynomials of a fixed degree which have
been obtained from other irreducible polynomials by applying a
fixed quadratic transformation. In this paper we show that in fact this is the case and 
for any fixed quadratic transformation, we determine the
number of irreducible polynomials of degree $n$ over a
finite field whose transformation is also irreducible. Our result
generalizes Carlitz's result---Carlitz's formula for the number of srim
polynomials can be recovered from our result. 

This paper is organized as follows. In Section~\ref{Main-result}, we state main result 
of the paper. In Section~\ref{Prel.}, we gather 
some results which will be used in the proof of main result. Section~\ref{proof}, 
contains the proof of the main result and finally we conclude by Section~\ref{com} 
which contains some remarks and comments.

In this paper $\cT(m,q)$ denotes the set of elements from $\fqm$ which are not in any
proper subfield of $\fqm$. 
\section{Main result}\label{Main-result}
Throughout this section and Section~\ref{proof} we assume that 
$g(x)=a_1x^2+b_1x+c_1$ and $h(x)=a_2x^2+b_2x+c_2$.  
\begin{theorem}
Let $q$ be a prime power, and let $g(x)$ and $h(x)\in \fq[x]$
be relatively prime polynomials with $\max(\deg(g),\deg(h))=2$.
Also let $I_{(n,g,h)}$ be the set of monic irreducible polynomials
$f(x)$ of degree $n > 1$ over $\fq$ whose quadratic transformation by $g(x)$
and $h(x)$, i.e. $p(x)=h(x)^nf\left(\frac{g(x)}{h(x)}\right)$, is irreducible
over $\fq$. Then
\[
\#I_{(n,g,h)} = \left\{
    \begin{array}{clclcl}
    0,& \;\;\;\mbox{if\;} b_1=b_2=0\;\mbox{and}\;q=2^l,\\\\
    \frac{1}{2n}(q^n-1), & \;\;\;\mbox{if\;q\;is\;odd\;and\;} n=2^m, m\ge 1, \\\\
    \frac{1}{2n}\displaystyle\sum_{\stackrel{d\mid n,}{d\;\mbox{{\tiny odd}}}}\mu(d) q^{n/d}, &\;\;\;\mbox{otherwise}.\\
    \end{array}
    \right.
\]
\end{theorem}
Notice that Carlitz' result follows from the above theorem for $g(x)=x^2+1$ and 
$h(x)=x$.

\section{Preliminaries}\label{Prel.}
In this section we gather some results which will be used in the rest
of the paper to prove the main result of the paper.

\subsection{Polynomial transformation of irreducible polynomials}
The following lemma is known as Capelli's lemma and can be found in ~\cite{Cohen-1} 
too.
\begin{lemma}\label{cohen}
Let $f(x)$ be a degree $n$ irreducible polynomial over $\mathbb{F}_q$, and
let $g(x), h(x)\in \mathbb{F}_q[x]$. Then $p(x)=h(x)^nf\left(g(x)/h(x)\right)$ is irreducible over 
$\mathbb{F}_q$ if and only if for any root $\alpha$ of $f(x)$ in 
$\mathbb{F}_{q^n}$, $g(x)-\alpha h(x)$ is an irreducible polynomial over 
$\mathbb{F}_{q^n}$.
\end{lemma}

\subsection{Resultant and discriminant}
We begin this section by recalling  the
{\it Resultant} of polynomials over a field. For a more detailed treatment
 see \cite[Ch. 1, pp. 35-37]{Harald-Rudolf}.

Let $F(x)$ and $G(x) \in K[x]$ and suppose $F(x)= a \prod_{i=0}^{s-1}(x-x_i)$
and $G(x) = b \prod_{j=0}^{t-1}(x-y_j)$, where $a,b \in K$ and
$x_0,x_1,\ldots,x_{s-1}$, $y_0,y_1,\ldots,y_{t-1}$ are in some
extension of $K$. Then the {\it Resultant}, $\Res(F,G)$, of $F(x)$
and $G(x)$ is

\begin{equation}\label{roots}
\Res(F,G) = (-1)^{st} b^s \prod_{j=0}^{t-1} F(y_j)
           = a^t \prod_{i=0}^{s-1} G(x_i) \; .
\end{equation}

Notice that $\Res(F,G)=0$ if and only if $F(x)$ and $G(x)$ have a common root in some extension
of $K$. Thus if $\Res(F,G)\neq 0$, then $F(x)$ and $G(x)$ are coprime. The following lemma 
which is probably already somewhere in the literature will be needed later. It follows from
direct calculations.
\begin{lemma}\label{Disc-xy}
Suppose $g(x)=a_1x^2+b_1x+c_1$, $h(x)=a_2x^2+b_2x+c_2$ and let $y$ be an indeterminate. Then
$$
\Disc_y(\Disc_x(g(x)-yh(x)))=16\Res(g,h),
$$
where $\Disc_x$ and $\Disc_y$ indicate that discriminant is taken with respect to variables 
$x$ and $y$, respectively.  
\end{lemma}
\subsection{Hurwitz genus formula}
\begin{lemma}\cite[Theorem~5.9]{Silverman}\label{Hurwitz} Let $\PP({\overline{\fq}})$ denote
 the one dimensional projective space over
 $\overline{\fq}$ (algebraic closure of $\fq$)
, and let $\phi:\PP({\overline{\fq}})\longrightarrow \PP({\overline{\fq}})$ be a 
a non-constant separable map. Then
\[
2\deg\phi-2\ge \sum_{P\in \PP}(e_{\phi}(P)-1),
\]
where $e_{\phi}(P)$ is the ramification index of $\phi$ at $P$. Equality holds if
and only if $e_\phi(P)$ is not divisible by the characteristic of $\fq$ for all
$P\in\PP({\overline{\fq}})$.
 \end{lemma}

\section{Proof of the main theorem}\label{proof}
\begin{proof}
Let $f(x)$ be a monic irreducible polynomial of degree $n$ over $\fq$. 
Using Theorem~\ref{cohen}, $p(x)$ is irreducible over $\fq$ if and only if for any root 
$\alpha$ of $f(x)$ in $\fqn$, $g(x)-\alpha h(x)$ is irreducible over $\fqn$. Thus in order
to compute the number of irreducible polynomials $p(x)=h(x)^nf\left(\frac{g(x)}{h(x)}\right)$ 
over $\fq$ 
it suffices to compute the number of elements $\beta\in \cT(n,q)$ for 
which 
$$r_\beta(x)=g(x)-\beta h(x)$$
is irreducible in $\fqn[x]$ and divide the result by $n$.

Now for $m=1,2,3,\ldots$, let 
\begin{equation}
U(m,q)=\left\{\beta: \beta\in\cT(m,q),\; r_\beta(x)\;\mbox{irreducible and quadratic over} 
\;\fqn\right\}.
\end{equation}

Notice that
$$
\#I_{(n,g,h)}= \frac{1}{n}\#U(n,q).
$$
It turns out that we first need to compute the number of elements of some auxiliary sets.
Thus let
\begin{equation}
\overline{U}(n,q)=\left\{\beta: \beta\in\fqn,\; r_\beta(x)\;\mbox{is irreducible and quadratic 
over} 
\;\fqn\right\},
\end{equation}
and
\begin{equation}
V(n,q)=\left\{\beta: \beta\in\fqn \;\mbox{and}\; \exists \gamma\in \fqn\; \mbox{s.t.}\; 
g(\gamma)=\beta h(\gamma)\right\}.
\end{equation}

If for $\beta\in\fqn$, $r_\beta(x)$ is not a constant polynomial over $\fqn$ , then it is either irreducible over $\fqn$ or it has exactly two roots
in $\fqn$ and gets factored as a product of two linear polynomials over $\fqn$. Thus
if we let $c$ be the number of $\beta\in\fqn$ for which $r_\beta(x)$ is a constant polynomial 
over $\fqn$, then 
\begin{equation}
\#V(n,q)=q^n-c-\#\overline{U}(n,q).
\end{equation}

In order to compute $\#V(n,q)$ we do a double counting as follows.
Suppose that 
\begin{equation}
W(n,q)=\big\{(\gamma,\beta): \gamma, \beta\in \fqn ; g(\gamma)=\beta h(\gamma)\big\}.
\end{equation}

For every $\gamma\in\fqn$, there is
a unique $\beta\in\fqn$ so that $g(\gamma)=\beta h(\gamma)$ unless $\gamma$ is a root of
$h(x)=0$. Now let the number of roots of $h(x)=0$ be $a$. Since $\deg(h(x))\le 2$, we have 
$a\le 2$. Thus 
\begin{equation}\label{count:W}
\#W(n,q)=q^n-a;\; a\le 2.
\end{equation}

On the other hand for every $\beta\in V(n,q)$ there are either
one or two $\gamma\in\fqn$ such that $g(\gamma)=\beta h(\gamma)$. In order to compute
 $\#V(n,q)$ we need to know how many elements 
in $V(n,q)$ have exactly one preimage and how many have exactly two preimages. We deal with
the fields of odd and even characteristic separately:
\begin{itemize} 
\item Fields of even characteristic: we have two cases
\begin{itemize}
\item $b_1=b_2=0$ : 
for every $\beta\in V(n,q)$ there is exactly one 
$\gamma\in\fqn$ such that $g(\gamma)=\beta h(\gamma)$. In fact in this case $p(x)$ is 
always square of a polynomial over $\fq$ and hence
the transformation by $g(x)$ and $h(x)$ does not result in any irreducible 
polynomial. This proves one of the cases of the main theorem. So in the rest of the 
paper we prove the remaining cases.
\item either $b_1$ or $b_2$ is nonzero: for every $\beta\in V(n,q)$ 
the equation   
$$
r_\beta(x)=g(x)-\beta h(x)=(a_1-\beta a_2)x^2+(b_1-\beta b_2)x+(c_1-\beta c_2)=0
$$
has exactly
two solutions in $\fqn$ unless either $r_\beta(x)$ is a linear polynomial or $b_1-\beta b_2=0$.
Each case can happen for at most one value of $\beta$.
\end{itemize}
\item Fields of odd characteristic:
for every $\beta\in V(n,q)$ there are exactly two $\gamma$ such that
  $g(\gamma)=\beta h(\gamma)$ unless either $r_\beta(x)$ is a linear polynomial or 
the discriminant of the equation
$$
r_\beta(x)=(a_1-\beta a_2)x^2+(b_1-\beta b_2)x+(c_1-\beta c_2)=0
$$
which is 
\begin{equation}\label{Disc:beta}
w(\beta)=\Disc(h)\beta^2+(4a_1c_2+4c_1a_2-2b_1b_2)\beta+\Disc(g).
\end{equation}
is zero. The former case can happen for at most one value of $\beta$. Now in the latter 
case we claim that Equation~\eqref{Disc:beta} can be zero for at most two values of $\beta$. 
It suffices to show that $w(\beta)$ is not identically zero. Suppose that $w(\beta)$ is 
identically zero. Then all the coefficients of $w(\beta)$ are zero and hence 
$\Disc(w(\beta))=0$. But using Lemma~\ref{Disc-xy} we have
\[
\Disc(w(\beta))=16\Res(g,h).
\]
This is a contradiction since we have assumed that $g(x)$ and $h(x)$ are relatively 
prime and do not have a common root.
\end{itemize}
Now if we let $d$ be the number of $\beta$ for which $r_\beta(x)$ is a linear polynomial
and $b$ be the number of $\beta$ for which $r_\beta(x)$ is a quadratic polynomial and $\beta$ has
one preimage, then $d\le 1$ and 
above arguments show that $b\le 2$. From this fact and Equation~\eqref{count:W} 
we deduce that 
$$
\#V(n,q)=\frac{\#W(n,q)+b+d}{2}=\frac{q^n-a+b+d}{2}
$$
or equivalently
\begin{equation}\label{count:U}
\#\overline{U}(n,q)=\frac{q^n+a-b-d-2c}{2}=\frac{q^n+a-b-2c-d}{2}.
\end{equation}

Having computed $\#\overline{U}(n,q)$ we can compute $\#U(n,q)$. We claim that
\begin{equation}\label{sum:U}
\#\overline{U}(n,q)\;=\;\sum_{\stackrel{d\mid n}{d\;{odd}}}\#U(n/d,q).
\end{equation}
Notice that if $\beta\in U(m,q)$ and $m\mid n$, then $\beta\in \overline{U}(n,q)$ if and
only if $n/m$ is an odd number since otherwise $\fqn$ contains a quadratic extension of
$\fqm$ and hence $r_\beta(x)$ is not irreducible over $\fqn$ any more. This 
proves~\eqref{sum:U}. 

Using~\eqref{count:U} and applying Mobius inversion we get
\begin{eqnarray}\label{count:mu}
\nonumber\#U(n,q)&=&\frac{1}{2}\sum_{\stackrel{d\mid n}{d\;{odd}}}\mu(d)(q^{n/d}+a-b-2c-d)\\
&=&\frac{1}{2}\sum_{\stackrel{d\mid n}{d\;{odd}}}\mu(d)q^{n/d}+
\frac{a-b-2c-d}{2}\sum_{\stackrel{d\mid n}{d\;{odd}}}\mu(d).
\end{eqnarray}
Now if $n$ is not a power of two and hence has at least two odd positive divisors, 
then we have
\[
\sum_{\stackrel{d\mid n}{d\;{odd}}}\mu(d)=0,
\]
and thus
\[
\#U(n,q)=\frac{1}{2}\sum_{\stackrel{d\mid n}{d\;{odd}}}\mu(d)q^{n/d}.
\]
This finishes the proof for the case of $n$ not being a power of two.
On the other hand if $n=2^m$ for some $m\ge 1$, then
\[
\#U(n,q)=\frac{1}{2}(q^n+a-b-2c-d).
\]
In the sequel, we show that $a-b-2c-d=-1$ if $\fq$ is of odd characteristic and
$a-b-2c-d=0$ if $\fq$ is of even characteristic. In order to prove this we can use  
elementary arguments and consider many cases but here we use Hurwitz genus formula
, Theorem~\ref{Hurwitz}, to give a shorter proof. 

Now let $\Phi$ be a map defined over one dimensional projective space over 
$\overline{\fq}$ as follows:
\begin{eqnarray}
\Phi&:&\PP({\overline{\fq}})\longrightarrow \PP({\overline{\fq}})\nonumber\\
\Phi([X:Y])&=&[a_1X^2+b_1XY+c_1Y^2:a_2X^2+b_2XY+c_2Y^2]\nonumber.
\end{eqnarray}

Over the fields of odd characteristic, $\Phi$ is obviously separable and 
non-constant and furthermore since over fields of even characteristic  
we have assumed that either $b_1$ or $b_2$ is nonzero
, it is separable over the fields of even characteristic, too. Thus one can apply 
Theorem~\ref{Hurwitz} to the map $\Phi$ and conclude that over fields of odd characteristic
$\Phi$ has two ramification points and over fields of even characteristic it has just 
one ramification point.

Now in order to avoid confusion in the rest of the proof, let $\infty_1$ and $\infty_2$ 
denote the points at infinity at the domain and range of the map $\Phi$, respectively.
In this setting, as we are assuming that $n=2^m$ for some $m\ge 1$, all the $x$ and $\beta$
related to $a,b,c,d$ are in $\fqtwo$ and thus $a$ is the number 
of finite preimages of $\infty_2$, $b$ is the number of finite ramification points with 
finite image, $c$ is one if $\infty_1$ is a ramification point and it image is finite and
zero otherwise, and $d$ is the number of finite points which have two
 preimages one of them being $\infty_1$. 

Suppose that $\fq$ is of odd characteristic. If $\infty_2$ is a branch point (its preimage
is a ramification point), then its preimage $\Phi^{-1}(\infty_2)$ is either $\infty_1$ or finite.  
If $\Phi^{-1}(\infty_2)$ is $\infty_1$, then $a=0$, $b=1$ as we can have one more ramification
point, $c=0$ as there is no finite branch point having $\infty_1$ as preimage and $d=0$
as there is no finite non-branch point having $\infty_1$ as preimage. If $\Phi^{-1}(\infty_2)$ 
is finite, then $a=1$ and there is one more ramification point. If $\infty_1$ is ramified,
then it is mapped to a finite point and hence $b=0$, $c=1$ and $d=0$ . If $\infty_1$ is 
unramified, then its image is finite and hence $b=1$ as there should be two ramified points,
$c=0$ and $d=1$. 

If $\infty_2$ is not a branch point, then either it has two finite preimages or it has 
one finite preimage. If it has two finite preimage, then $a=2$, either $b=1$, $c=1$ and 
$d=0$ if $\infty_1$ is a ramification point or $b=2$, $c=0$ and $d=1$ if $\infty_1$ is not 
a ramification point. If $\infty_2$ has one finite preimage, then $a=1$, $b=2$, $c=0$ and
$d=0$. 

We see that in all the cases if the characteristic of $\fq$ is an odd number, then
$a-b-2c-d=-1$.

Similar arguments show that in the case of fields of even characteristic $a-b-2c-d=0$. 
This finishes the proof as if $n=2^m$ and $q$ is a power of two then
$$
\frac{1}{2}\sum_{\stackrel{d\mid n}{d\;{odd}}}\mu(d)q^{n/d}=\frac{1}{2}q^n.
$$ 
\end{proof}

\section{Comments}\label{com}
\subsection{Alternative proof approach}
Since $r_\beta(x)$ is irreducible over fields of odd characteristic if and only if 
its discriminant is a non-square in $\fqn$,
another approach that can be used to prove the main result for fields of odd characteristic
is to see for how many $\beta$ discriminant of $r_\beta(x)$ is a non-square and for 
how many $\beta$ it is a square in $\fqn$. This can be done using the following well-known
lemma. The following lemma implies that a quadratic polynomial over a finite field of odd 
characteristic is square almost as many times as it is a non-square. 
\begin{lemma}\cite[Theorem~5.48]{Harald-Rudolf}\label{quad-sum}
Let $q$ be an odd prime power, and let
$f(x)=a_2x^2+a_1x+a_0\in\fq[x]$ where $a_2\neq 0$. Let $\eta$ be
the quadratic character of $\fqn$. If $\Disc(f)\neq 0$, then
$\sum_{c\in\fqn}\eta(f(c))=-\eta(a_2)$.
\end{lemma}

\subsection{Palindromic primes}
It is well known that the number $I(2n,q)$ of irreducible polynomials of degree $2n$ 
over $\fq$ is
\[
I(2n,q)=\frac{1}{2n}\sum_{\stackrel{d\mid 2n}{}}\mu(d)q^{2n/d},
\]
and the probability that a random monic polynomial of degree $2n$ is irreducible over $\fq$
is roughly $\frac{1}{2n}$. On the other hand number of polynomials of degree $2n$ obtained from a 
fixed quadratic transformation is $q^n$. Thus Carlitz' result and  our result implies that 
the number of irreducible polynomials among the polynomials obtained by a quadratic 
transformation is roughly what one would expect. In~\cite{billbanks}, it was shown that 
\[
\frac{\mbox{Number of palindromic primes} \le x\; {\mbox{written in base} }\;g}
{\mbox{Number of palindromic numbers} \le x\; {\mbox{written in base} }\;g}=O\left(\frac{\log\log\log x}
{\log\log x}\right)
\]
where the implied constant depends only on the base $g$ and conjectured that
\[
\frac{\mbox{Number of palindromic primes} \le x\; {\mbox{written in base} }\;g}
{\mbox{Number of palindromic numbers} \le x\; {\mbox{written in base} }\;g}\sim C\frac{1}
{\log x}
\] 
or roughly speaking palindromic numbers with respect to primality behave like random integers.
Carlitz' result can be viewed as an affirmative answer to their conjecture in the finite field
setting. Now we wonder what the analogue of our result is for the integer numbers and if 
one can establish results similar to those of~\cite{billbanks}?

\section{Acknowledgments}
The author would like to thank Vijaykumar Singh for helpful discussion during the 
preparation of this paper and Stephen Cohen and Igor Shparlinski for comments on an earlier
draft of it.

Research of the author is supported by the Claude Shannon Institute,
 Science Foundation Ireland Grant 06/MI/006.
\bibliographystyle{plain}
\bibliography{FF}

\end{document}